\documentclass[pdflatex,sn-mathphys-num]{sn-jnl}

\usepackage{graphicx}%
\usepackage{multirow}%
\usepackage{amsmath,amssymb,amsfonts}%
\usepackage{amsthm}%
\usepackage{mathrsfs}%
\usepackage[title]{appendix}%
\usepackage{xcolor}%
\usepackage{textcomp}%
\usepackage{manyfoot}%
\usepackage{booktabs}%
\usepackage{algorithm}%
\usepackage{algorithmicx}%
\usepackage{algpseudocode}%
\usepackage{listings}%

\usepackage{multicol,multirow}
\usepackage{rotating}

\usepackage{ifpdf}
\usepackage[T1]{fontenc}
\usepackage{newtxtext}
\usepackage{newtxmath}
\usepackage{textcomp}

\usepackage{comment}
\usepackage[numbers]{natbib}
\usepackage{mathtools}
\usepackage{yfonts}
\usepackage{tikz-cd}

\usepackage{dirtytalk}

\theoremstyle{thmstyleone}%
\newtheorem{theorem}{Theorem}[section]
\theoremstyle{thmstyletwo}%

\theoremstyle{thmstylethree}%
\newtheorem{definition}{Definition}%
\newtheorem{conjecture}{Conjecture}

\newtheorem*{theoremA*}{Theorem A}
\newtheorem*{theoremB*}{Theorem B}
\newtheorem*{theoremC*}{Theorem C}
\newtheorem*{theoremD*}{Theorem D}

\newtheorem*{theorem*}{Theorem}

\DeclareFontFamily{U}{rcjhbltx}{}
\DeclareFontShape{U}{rcjhbltx}{m}{n}{<->rcjhbltx}{}
\DeclareSymbolFont{hebrewletters}{U}{rcjhbltx}{m}{n}
\DeclareMathSymbol{\shin}{\mathord}{hebrewletters}{152}
\DeclareMathSymbol{\samech}{\mathord}{hebrewletters}{115}

\newcommand{\pt}[2]{\left\langle\!\left\langle #1,#2 \right\rangle\!\right\rangle}

\begin{document}

\title[Article Title]{Shimura Reciprocity for Real Quadratic Fields}


\author*[1]{\fnm{Joseph} \sur{DiCapua}}\email{jdicap1@lsu.edu}

\affil*[1]{\orgdiv{ORCID 0009-0004-5894-9523, Department of Mathematics}, \orgname{Louisiana State University}, \orgaddress{\street{175 Field House Drive}, \city{Baton Rouge}, \postcode{70802}, \state{LA}, \country{United States}}}


\abstract{The purpose of this note is to give conditional explicit formulas for the Galois action on real multiplication values of the Shintani--Faddeev modular cocycle. We assume the order 1 abelian Stark conjecture for real quadratic fields. We conclude this note by presenting our own conjecture about the Galois action. Numerical evidence for this conjecture is provided by the algorithm described in \hyperref[sec4]{Section 4} and the attached code.}




\maketitle 

\section{Introduction}
A SIC (or symmetric informationally complete positive operator valued measure) can be described as a collection of $d^2$ equiangular lines all going through the origin in $\mathbb{C}^d$. Equivalently such a construction is a complex equiangular tight frame for which the upper bound of $d^2$ vectors in dimension $d$ is achieved \cite{ApplebyFlammiaKopp, DelsarteGoethalsSeidel}. The question of whether SICs exist in all finite dimensions is an important problem in quantum theory known as Zauner's conjecture, and it arises as a kind of engineering problem in classical signal processing \cite{ApplebyBengtsson2019}. Recently there has been significant interest and progress on this conjecture in \cite{ApplebyFlammiaKopp} and \cite{RadchenkoWheeler1}.
\\
\\
\noindent
On the other hand the Stark conjectures \cite{Stark1, Stark2, Stark3, Stark4, Tate1} concern the properties of special values of derivatives of zeta functions in algebraic number theory. They are closely related to Hilbert's twelfth problem. In \cite{ApplebyFlammiaKopp}, the authors show Zauner's conjecture follows from Stark conjectures together with a special function identity called the \textit{Twisted Convolution Conjecture}. This special function identity was recently proved in \cite{RadchenkoWheeler1} and \cite{ApplebyFlammiaKopp2}.
\\
\\
\noindent
Conjecture 1.50 of \cite{ApplebyFlammiaKopp} states that the ghost SIC overlaps generate certain ray class fields of a quadratic extension $K$. Our paper makes heavy use of Theorem 2.19 from the same paper in order to get a conditional explicit description of the action of the correct ray class group of $K$ on ghost SIC overlaps (Stark units for $K$). We also present an algorithm based on this description which can be used to compute the Galois action on a ghost SIC arising from an admissible tuple as in \cite{ApplebyFlammiaKopp} for any dimension $d$.
\\
\\
\noindent
\textbf{Acknowledgments:} I would like to thank Gene Kopp for introducing me to this problem.

\section{Background and notation}

Let $K$ be a fixed real quadratic field. First we recall the definition of an admissible tuple following \cite{ApplebyFlammiaKopp}. In doing this we recall the notions of dimension grid and rank grid associated to $K$. \begin{definition}
We define the fundamental totally positive unit $\epsilon$ of $K$ to be the smallest positive norm unit of $K$ which is greater than $1$.
\end{definition}

\noindent
We also define the sequence of conductors following \cite{ApplebyFlammiaKopp}. \begin{definition}
Let $\Delta_0$ be the discriminant of $K$. For positive integers $j$ define

\[
f_j = \frac{\epsilon^j - \epsilon^{-j}}{\sqrt{\Delta_0}}
\]
\end{definition}

\noindent
Note each $f_j$ is a positive integer. We now define three sequences or \say{grids} associated to $K$. \begin{definition}
For all positive integers $j,m$ define $r_{j,m} = f_{jm}/f_j$. We also define $d_{j,m} = r_{j,m+1} + r_{j,m}$ and $d_j = d_{j,1}$. The two-index sequence $d_{j,m}$ is the dimension grid associated to $K$ and $r_{j,m}$ is the rank grid associated to $K$. The sequence $d_j$ is the dimension tower associated to $K$. 
\end{definition}

\noindent
Now for positive integers $j,m$ we say the triple $(K,j,m)$ is an admissible triple. At this point we also need to give the definition of an admissible pair
\begin{definition}
    A pair of integers $(d,r)$ is called admissible if there exists an $n > 4$ such that 
    \begin{itemize}
        \item $ 0 < r < (d-1)/2 $ and
        \item $ nr(d-r) = d^2-1 $
    \end{itemize}
    For each admissible pair we say that $L = \mathbb{Q}(\sqrt{n(n-4)})$ is the real quadratic field associated to the pair.
\end{definition}

\begin{theorem}
There is a bijection between admissible triples and admissible pairs. For each admissible triple $(K,j,m)$, the pair $(d_{j,m},r_{j,m})$ is admissible. Conversely, for each admissible pair, $(d,r)$ there is a unique admissible triple $(K,j,m)$ such that $(d_{j,m},r_{j,m}) = (d,r)$.  
\end{theorem}

\noindent
The above is Theorem 1.25 from \cite{ApplebyFlammiaKopp}. We also have
\begin{definition}
    Let $Q$ be an integral binary quadratic form, and let $(d,r) \sim (K,j,m)$ be admissible tuples. We say that $(d,r,Q)\sim (K,j,m,Q)$ are admissible if the fundamental discriminant of $Q$ is the discriminant of $K$ and the conductor of $Q$ is a divisor of $f_j$.
\end{definition}
\noindent
We recall another theorem from \cite{ApplebyFlammiaKopp} which is also Theorem 1.1 from \cite{KoppLagarias} which allows us to discuss ray class fields of $K$.

\begin{theorem}
Let $K$ be a number field, and let $(\mathcal{O};\mathfrak{m},\Sigma)$ be a level datum for $K$. Then there exists a unique abelian Galois extension $H^{\mathcal{O}}_{\mathfrak{m}, \Sigma}/K$ with the property that a prime ideal $\mathfrak{p}$ of $\mathcal{O}_K$ that is coprime to the quotient ideal $(\mathfrak{m} : \mathcal{O}_K)$ splits completely in $H^\mathcal{O}_{\mathfrak{m},\Sigma}/K$ if and only if $\mathfrak{p} \cap \mathcal{O} = \pi \mathcal{O}$, a principal prime $\mathcal{O}$-ideal having $\pi \in \mathcal{O}$ with $\pi \equiv 1 (\text{mod } \mathfrak{m})$ and $\sigma(\pi) > 0$ for all $\sigma \in \Sigma$.

\noindent
Additionally, these fields have the following properties:
\begin{itemize}
\item $H^{\mathcal{O}_K}_{\mathfrak{m}\mathcal{O}_K,\Sigma} \subseteq H^\mathcal{O}_{\mathfrak{m}, \Sigma} \subseteq H^{\mathcal{O}_K}_{(\mathfrak{m}:\mathcal{O}_K),\Sigma}$.
\item There is a canonical isomorphism $\text{Art}_\mathcal{O} : \text{Cl}_{\mathfrak{m},\Sigma}(\mathcal{O}) \rightarrow \text{Gal}(H^\mathcal{O}_{\mathfrak{m},\Sigma}/K)$.
\end{itemize}
\end{theorem}

\noindent
We also need the definition of the \textit{flat imprimitive ray class monoid} from \cite{ApplebyFlammiaKopp}. This finite commutative monoid contains the ray class group as a submonoid, and it will be needed to define generalized zeta values that give rise to ghost SIC overlaps. 
\begin{definition}
    The \textit{flat imprimitive ray class monoid} is 
    \[
    \overline{\text{Clm}}^\flat_{\mathfrak{m},\Sigma}(\mathcal{O}) = \frac{\overline{J}^\flat_\mathfrak{m}(\mathcal{O})}{\sim_{\mathfrak{m},\Sigma}},
    \]
    where 
    \[
    \overline{J}^\flat_\mathfrak{m}(\mathcal{O}) = \{ \mathfrak{a} \in J^*_\mathcal{O}(\mathcal{O}) \mid \mathfrak{a}\mathcal{O}[S_\mathfrak{m}^{-1}] \subseteq \mathcal{O}[S_\mathfrak{m}^{-1}]\}
    \]
    with
    \[
    S_\mathfrak{m} = \{ \alpha \in \mathcal{O} \mid \alpha\mathcal{O} + \mathfrak{m} = \mathcal{O} \},
    \]
    and the equivalence relation $\sim_{\mathfrak{m},\Sigma}$ is defined by
    \begin{equation}
  \mathfrak{a} \sim_{\mathfrak{m},\Sigma} \mathfrak{b}
  \iff
  \left\{
  \begin{array}{l}
    \text{there exists } \mathfrak{c} \in \overline{J}^\flat_{\mathfrak{m}}(\mathcal{O})
    \text{ and } \alpha,\beta \in \mathcal{O}[S_{\mathfrak{m}}^{-1}] \text{ such that} \\[2pt]
    \quad \mathfrak{a} = \alpha\mathfrak{c}, \ \mathfrak{b} = \beta\mathfrak{c}, \\[2pt]
    \quad \alpha-\beta \in \mathfrak{m}\mathcal{O}[S_{\mathfrak{m}}^{-1}], \text{ and} \\[2pt]
    \quad \operatorname{sgn}(\sigma(\alpha)) = \operatorname{sgn}(\sigma(\beta))
      \text{ for all } \sigma \in \Sigma.
  \end{array}
  \right.
\end{equation}
In the above, $J_\mathfrak{m}^*(\mathcal{O}) = \{ \text{invertible fractional ideals of }\mathcal{O}\text{ coprime to }\mathfrak{m} \}$.
\end{definition}

\noindent
We need several more definitions from \cite{ApplebyFlammiaKopp} in order to discuss ghost overlaps of potential ghost SICs.

\begin{definition}
For all $M = \begin{pmatrix}
\alpha && \beta\\
\gamma && \delta
\end{pmatrix} \in \text{SL}_2(\mathbb{Z})$ the Rademacher class invariant is given by 
\begin{itemize}
    \item $\Psi(M) = \frac{\text{Tr}(M)}{\gamma} - 3\text{sgn}(\gamma\text{Tr}(M)) - 12\text{sgn}(\gamma)\mathfrak{s}(\alpha,\gamma)$ if $\gamma \neq 0$ and
    \item $\Psi(M) = \beta/\delta$ if $\gamma = 0$
\end{itemize}
\noindent
where $\mathfrak{s}(a,b)$ with $a,b \in \mathbb{Z}$ and $b \neq 0$ is the Dedekind sum
\[
\mathfrak{s}(a,b) = \sum_{n=1}^{|b|-1} \left(\left(\frac{n}{b}\right)\right)\left(\left(\frac{na}{b}\right)\right)
\]
\noindent
with $((x)) = 0$ if $x \in \mathbb{Z}$ and $((x)) = x - \lfloor x \rfloor - \frac{1}{2}$ otherwise. Here we also adopt the convention that $\text{sgn}(0) = 0$.
\end{definition}

\begin{definition}
Let $t = (d,r,Q) \sim (K,j,m,Q)$ be an admissible tuple, and let $\mathbf{p} \in \mathbb{Z}^2$. The SF phase, denoted $\phi_\mathbf{p}(t)$, is
\[
\phi_\mathbf{p}(t) = (-1)^{s_d(\mathbf{p})}e^{-\frac{\pi i}{12} \Psi(A_t)} \xi_d^{-\frac{f_{jm}}{f}Q(\mathbf{p})}
\]
\noindent
where $s_d(\mathbf{p}) = d+ (1+d)(1+p_1)(1+p_2)$, $f$ is the conductor of $Q$, and $\xi_d$ is a fixed primitive $\overline{d}$th root of unity where $\overline{d} = d$ if $d$ odd and $\overline{d} = 2d$ if $d$ even.
\end{definition}

\noindent
We refer to definition 1.14 of \cite{ApplebyFlammiaKopp} for the definition of different variant $q$-Pochhammer symbols.

\begin{definition}
The \textit{finite variant $q$-Pochhammer symbol} is defined by 

\[
\varpi_n(z,\tau) =
\begin{cases}
\prod_{j=0}^{n-1}\left(1-e^{2\pi i (z+j\tau)}\right) & \text{if } n > 0, \\
1  & \text{if } n = 0,\\
\prod_{j=n}^{-1}\left(1-e^{2\pi i (z+j\tau)}\right)^{-1} & \text{if } n < 0
\end{cases}
\]
for $n\in \mathbb{Z}$ and $z,\tau \in \mathbb{C}$. The (infinite) \textit{variant $q$-Pochhammer symbol} is defined By
\[
\varpi(z,\tau) = \prod_{j=0}^\infty\left(1-e^{2\pi i (z+j\tau)}\right)
\]
for $z\in\mathbb{C}$ and $\tau \in \mathbb{H}$.
\end{definition}

\noindent
For $\tau \in \mathbb{C}$ and $M = \begin{pmatrix} \alpha && \beta \\ \gamma && \delta \end{pmatrix} \in \text{GL}_2(\mathbb{Z})$, define the \textit{fractional linear transform} $M \cdot \tau$ and denote the denominator $j_M(\tau)$ respectively by
\[
M\cdot \tau = \frac{\alpha\tau + \beta}{\gamma\tau + \delta}
\]
and
\[
j_M(\tau) = \gamma\tau + \delta
\]
We say $\tau$ is a fixed point of $M$ if $M\cdot\tau = \tau$.

\begin{definition}
For $M = \begin{pmatrix} \alpha && \beta \\ \gamma && \delta \end{pmatrix} \in \text{GL}_2(\mathbb{Z})$ define $\mathcal{D}_M$ to be the set $\mathbb{C} \setminus \{ \tau \in \mathbb{R} \mid \text{det}(M)j_M(\tau) \le 0 \}$.
\end{definition}

\begin{definition}
For $M \in \text{SL}_2(\mathbb{Z})$ the \textit{Shintani--Faddeev Jacobi cocycle} is a meromorphic function $\sigma_M$ on $\mathbb{C}\times \mathcal{D}_M$ whose restriction to $\mathbb{C}\times \mathbb{H}$ is given by 
\[
\sigma_M(z,\tau) = \frac{\varpi\left(\frac{z}{j_M(\tau)}, M\cdot \tau\right)}{\varpi(z,\tau)}
\]
See \cite{Kopp1} and \cite{ApplebyFlammiaKopp} for the continuation to $\mathbb{C}\times\mathcal{D}_M$.
\end{definition}

\noindent
For $\mathbf{p} = \begin{pmatrix} p_1 \\ p_2 \end{pmatrix}$ and $\mathbf{q} = \begin{pmatrix} q_1 \\ q_2 \end{pmatrix} \in \mathbb{C}^2$ define the nondegenerate \textit{symplectic form}
\[
\langle \mathbf{p}, \mathbf{q} \rangle = p_2q_1 - p_1q_2
\]
When the second argument is a complex number $\tau \in \mathbb{C}$ we also use the notation
\[
\pt{\mathbf{p}}{\tau} = \langle \mathbf{p}, \begin{pmatrix} \tau \\ 1 \end{pmatrix} \rangle = p_2\tau - p_1
\]

\begin{definition}
For $\mathbf{r} \in \mathbb{Q}^2$, let $\Gamma_\mathbf{r}$ be the subgroup of $\text{SL}_2(\mathbb{Z})$ consisting of matrices $A$ such that $(A-I)\mathbf{r} \in \mathbb{Z}^2$.
\end{definition}

\noindent
At this point we refer to definition 1.18 of \cite{ApplebyFlammiaKopp} for the definition of the Shintani--Faddeev modular cocycle.

\begin{definition}
If $\mathbf{r}\in\mathbb{Q}^2$ and $M \in \Gamma_\mathbf{r}$, then the Shintani--Faddeev (SF) modular cocycle is defined by
\[
\shin_M^\mathbf{r}(\tau) = \frac{\sigma_M(\pt{\mathbf{r}}{\tau}, \tau )}{\varpi_{((I-M)\mathbf{r})_2}\left(\frac{\pt{\mathbf{r}}{\tau} }{j_M(\tau)},M\cdot \tau\right)}
\]
for all $\tau \in \mathcal{D}_M$.
\end{definition}

\noindent
Still following \cite{ApplebyFlammiaKopp} let $Q = \langle a,b,c \rangle$ be an integral binary quadratic form; then $\Delta = -4\text{det}(Q) = b^2-4ac$ is its discriminant. The roots of $Q$ are $\rho_{Q,\pm} = \frac{-b \pm\sqrt{\Delta}}{2a}$.

\begin{definition}
Let $t = (d,r,Q) \sim (K,j,m,Q)$ be an admissible tuple, and let $\mathbf{p} \in \mathbb{Z}^2$. The corresponding candidate normalized ghost overlaps are defined by
\[
\tilde{\nu}_{\mathbf{p}}(t) = \phi_\mathbf{p}\shin_{A_t}^{d^{-1}\mathbf{p}}(\rho_t)
\]
\noindent
for all $\mathbf{p}$, where $\rho_t = \rho_{Q,+}$. The corresponding candidate ghost overlaps are defined by
\[
\tilde{\mu}_\mathbf{p}(t) = r \text{ if } \mathbf{p} \equiv 0 \text{ mod } d \text{ and } \tilde{\mu}_\mathbf{p}(t) = \sqrt{\frac{r(d-r)}{d^2-1}}\tilde{\nu}_\mathbf{p}(t) \text{ otherwise.}
\]
\end{definition}

\section{Galois action on real multiplication values of $\shin$}
\label{sec3}

\noindent
Let $K$ be a fixed real quadratic field coming from some admissible tuple $t = (d,r,Q) \sim (K,j,m,Q)$. In what follows, $\mathcal{O}_f$ is the order of $K$ with conductor $f$, where $f$ is the conductor of $Q$. We use $H^{\mathcal{O}_f}_{\overline{d} \infty_2}$ to denote the ray class field with level datum $(\mathcal{O}_f; \overline{d}\mathcal{O}_f, \{\infty_2\})$. For arbitrary $h \in G(H^{\mathcal{O}_f}_{\overline{d} \infty_2}/K)$ we would like to understand how $h$ acts on the ghost overlaps $\tilde{\mu}_\mathbf{p}(t)$. We will start with some special subgroup first to simplify things.
\\
\\
Now

$$\tilde{\mu}_\mathbf{p}(t) = \frac{1}{\sqrt{d_j+1}} \phi_{\mathbf{p}}(t)\shin^{d^{-1}\mathbf{p}}_A(\rho_t)$$

\noindent
so in order to understand how $h$ acts on $\tilde{\mu}_\mathbf{p}(t)$ it would be very helpful to see how $h$ acts on $\shin^{d^{-1}\mathbf{p}}_A(\rho_t)$. By Theorem 2.19 in \cite{ApplebyFlammiaKopp} and the remarks directly after this theorem, there is a class $\textfrak{U}$ in $\overline{\text{Clm}}^\flat_{\overline{d} \infty_2}(\mathcal{O}_f)$ that maps to the $\text{SL}_2(\mathbb{Z})$-orbit of $(d^{-1}\mathbf{p}, \rho_t)$ under the map $\Upsilon_{\overline{d}\mathcal{O}_f}$. See the paragraphs directly after Theorem 2.19 in \cite{ApplebyFlammiaKopp} and Theorem 3.14 of \cite{Kopp1} for descriptions of this map. By the same theorem we have

$$(\psi^{-2}\chi_{d^{-1}\mathbf{p}}^{-1})(A_t)\shin^{d^{-1}\mathbf{p}}_A(\rho_t)^2 = \exp(nZ'_{\overline{d} \infty_2}(0,\textfrak{U}))$$

\noindent
where $(\psi^{-2}\chi_{d^{-1}\mathbf{p}}^{-1})(A_t) = \phi_\mathbf{p}(t)^2$. See definition 2.4 of \cite{ApplebyFlammiaKopp} for the definition of the \textit{differenced ray class partial zeta function} $Z_{\mathfrak{m},\Sigma}(s,\mathfrak{a})$. We also have $h \in G(H^{\mathcal{O}_f}_{\overline{d}\infty_2}/K)$ so there exists some $\textfrak{B} \in \text{Cl}_{\overline{d} \infty_2}(\mathcal{O}_f)$ such that $\text{Art}(\textfrak{B}) = h$. Then by the \say{Monoid Stark Conjecture} (Conjecture 2.9 of \cite{ApplebyFlammiaKopp}), we have

$$h(\exp(n Z'_{\overline{d} \infty_2}(0,\textfrak{U}))) = \exp(n Z'_{\overline{d}\infty_2}(0,\textfrak{B}\textfrak{U}))$$

\noindent
First we work out the case where $\textfrak{B}$ is principal, and then we work out the general case. Let $(\textfrak{BU})_0$ be the class of $\textfrak{BU}$ in $\text{Cl}(\mathcal{O}_f)$. Let $\textfrak{U}_0$ be the class of $\textfrak{U}$ in $\text{Cl}(\mathcal{O}_f)$. Note that $(\textfrak{BU})_0 = \textfrak{U}_0$ because $\textfrak{B}$ is principal. We have that we can choose some $\textfrak{b}_1 \in \textfrak{U}_0^{-1}$ coprime to $\overline{d} \mathcal{O}_f$ and we can write $\textfrak{b}_1\overline{d}\mathcal{O}_f = \alpha(\rho_t \mathbb{Z} + \mathbb{Z})$ for some $\alpha \in K$ such that $\alpha$ is totally positive. We also have the following: we can choose some $\textfrak{b}_2 \in (\textfrak{BU})_0^{-1}$ coprime to $\overline{d}\mathcal{O}_f$ for which we can write $\textfrak{b}_2d\mathcal{O}_f = \alpha'(\rho_t\mathbb{Z} + \mathbb{Z})$ for some $\alpha' \in K$, $\alpha'$ totally positive. 
\\
\\
\noindent
Next there is some 

\begin{align}
    \mathbf{r} &= \begin{bmatrix}
        r_1 \\
        r_2
    \end{bmatrix} \in \mathbb{Q}^2
\end{align} 

\noindent
such that $(\alpha(r_2\rho_t - r_1))\textfrak{b}_1^{-1} \in \textfrak{U}$ and $r_2\rho' - r_1 > 0$. $(\mathbf{r},\rho_t)$ needs to be in the same $\text{SL}_2(\mathbb{Z})$-orbit as $(d^{-1}\mathbf{p}, \rho_t)$. Here $\rho'$ is the conjugate of $\rho$.
\\
\\
We can also choose some

\begin{align}
    \mathbf{s} &= \begin{bmatrix}
        s_1 \\
        s_2
    \end{bmatrix} \in \mathbb{Q}^2
\end{align} 

\noindent
such that $(\alpha'(s_2\rho_t - s_1))\textfrak{b}_2^{-1} \in \textfrak{BU}$.
\\
\\
Pick any $\textfrak{b} \in \textfrak{B}$. Then $\textfrak{b}(\alpha(r_2\rho_t - r_1))\textfrak{b}_1^{-1} \in \textfrak{BU}$. The idea now is to pick some $b \in \textfrak{b}$ such that

$$b(\alpha(r_2\rho_t - r_1))\textfrak{b}_1^{-1} \in \textfrak{BU}$$

\noindent
we are allowed to do this because $\textfrak{B}$ maps to the trivial class in $\text{Cl}(\mathcal{O}_f)$, hence $\textfrak{b}$ is principal.
Multiplication by $b$ in $K$ is given by some matrix $M_b$ and we have

$$b(\alpha(r_2\rho_t - r_1))\textfrak{b}_1^{-1} = \alpha(r_2'\rho_t - r_1')\textfrak{b}_1^{-1}$$

\noindent
where

\begin{align}
    \mathbf{r}' = \begin{bmatrix}
        r_1' \\
        r_2'
    \end{bmatrix} = M_{b}
    \begin{bmatrix}
        r_1 \\
        r_2
    \end{bmatrix}
\end{align} 

\noindent
One additional condition we need to check is $r'_2\rho_t' - r_1' > 0$. If this does not hold, just replace $b$ with $-b$. Then $\Upsilon_{\overline{d}\mathcal{O}_f}$ maps the class $\textfrak{BU}$ in $\overline{\text{Clm}}^\flat_{\overline{d} \infty_2}(\mathcal{O}_f)$ to the $\text{SL}_2(\mathbb{Z})$-orbit of $(\mathbf{r}', \rho_t)$ and we get

$$(\psi^{-2}\chi_{\mathbf{r}'}^{-1})(A_t)\shin^{\mathbf{r}'}_A(\rho_t)^2 = \exp(nZ'_{\overline{d} \infty_2}(0,\textfrak{BU}))$$

\noindent
We can conclude that

$$h(\phi_\mathbf{p}(t)^2\shin_A^{d^{-1}\mathbf{p}}(\rho_t)^2) = \phi_{M_b\mathbf{p}}(t)^2\shin_A^{d^{-1}M_b\mathbf{p}}(\rho_t)^2$$

\noindent
We get the following by taking square roots:

\begin{theorem}
Assume the order 1 abelian Stark conjecture for real quadratic fields. Then using the notation defined above we have
$$h(\phi_\mathbf{p}(t)\shin_A^{d^{-1}\mathbf{p}}(\rho_t)) = \pm\phi_{M_b\mathbf{p}}(t)\shin_A^{d^{-1}M_b\mathbf{p}}(\rho_t)$$
whenever $h$ corresponds to a principal ideal in $\text{Cl}_{\overline{d} \infty_2}(\mathcal{O}_f)$.
\end{theorem}

\noindent
See Conjecture 1 at the end of this note for how we believe the above sign is determined. 

\noindent
Now we take arbitrary $h \in G(H^{\mathcal{O}_f}_{d \infty_2}/K)$ and describe the action of $h$ on $\shin_A^{d^{-1}\mathbf{p}}(\rho_t)$. Using the same choice of $\textfrak{U} \in \overline{\text{Clm}}^\flat_{\overline{d} \infty_2}$ from before we still have

$$(\psi^{-2}\chi_{d^{-1}\mathbf{p}}^{-1})(A_t)\shin^{d^{-1}\mathbf{p}}_A(\rho_t)^2 = \exp(nZ'_{\overline{d} \infty_2}(0,\textfrak{U}))$$

\noindent
There exists some $\textfrak{B} \in \text{Cl}_{\overline{d} \infty_2}(\mathcal{O}_f)$ such that $\text{Art}(\textfrak{B}) = h$. Once again we have

$$h(\exp(n Z'_{\overline{d} \infty_2}(0,\textfrak{U}))) = \exp(n Z'_{\overline{d}\infty_2}(0,\textfrak{B}\textfrak{U}))$$

\noindent
At this point we fix a set of representatives $\{ \mathcal{O}_f = \mathbb{Z}\rho^{(1)} + \mathbb{Z}, \ldots ,\mathbb{Z}\rho^{(h_f)}+ \mathbb{Z}\}$ for $\text{Cl}(\mathcal{O}_f)$. For each of these representatives we also fix a class $\textfrak{B}_{\rho^{(i)}} \in \overline{\text{Clm}}^\flat_{\overline{d} \infty_2}(\mathcal{O}_f)$ such that the class of $\textfrak{B}_{\rho^{(i)}}^{-1}$ in $\text{Cl}(\mathcal{O}_f)$ is $[\mathbb{Z}\rho^{(i)}+\mathbb{Z}]$. Since we have already established how $h$ acts on $\shin_A^{d^{-1}\mathbf{p}}(\rho_t)$ if $h \in G(H^{\mathcal{O}_f}_{d \infty_2}/H^{\mathcal{O}_f})$ where $H^{\mathcal{O}_f}$ is the ring class field of $\mathcal{O}_f$, we may assume our $\textfrak{B} = \textfrak{B}_{\rho^{(i)}}$ for some $i$. We can also assume $\rho_t = \rho^{(i_t)}$.
\\
\\
Now we can pick some $\textfrak{b} \in \textfrak{U}_0^{-1}$ coprime to $\overline{d}\mathcal{O}_f$ and write $\textfrak{b}\overline{d}\mathcal{O}_f = \alpha(\mathbb{Z}\rho^{(i_t)} + \mathbb{Z})$ such that $\alpha$ is totally positive. We also have

\begin{align}
    \mathbf{r} &= \begin{bmatrix}
        r_1 \\
        r_2
    \end{bmatrix} \in \mathbb{Q}^2
\end{align} 

\noindent 
such that $\alpha(r_2\rho^{(i_t)} - r_1)\textfrak{b}^{-1} \in \textfrak{U}$. Let $\textfrak{B}_0$ be the class of $\textfrak{B} = \textfrak{B}_{\rho^{(i)}}$ in $\text{Cl}(\mathcal{O}_f)$, so $\textfrak{B}^{-1}_0 = [\mathbb{Z}\rho^{(i)} + \mathbb{Z}]$. Pick some $\textfrak{b}' \in \textfrak{B}_0^{-1}$ coprime to $\overline{d}\mathcal{O}_f$ and write $\textfrak{b}'\overline{d}\mathcal{O}_f = \alpha'(\mathbb{Z}\rho^{(i)} + \mathbb{Z})$ such that $\alpha'$ is totally positive. We have some

\begin{align}
    \mathbf{s} &= \begin{bmatrix}
        s_1 \\
        s_2
    \end{bmatrix} \in \mathbb{Q}^2
\end{align} 

\noindent
such that $\alpha'(s_2\rho^{(i)} - s_1)\textfrak{b}'^{-1} \in \textfrak{B}_{\rho^{(i)}}$. The class of $(\textfrak{B}_{\rho^{(i)}}\textfrak{U})^{-1}$ in $\text{Cl}(\mathcal{O}_f)$ will be the product $[\mathbb{Z}\rho^{(i_t)} + \mathbb{Z}][\mathbb{Z}\rho^{(i)} + \mathbb{Z}] = [\mathbb{Z}\rho^{(j)} + \mathbb{Z}]$. Then the product $\alpha\alpha'(r_2\rho^{(i_t)} - r_1)(s_2\rho^{(i)} - s_1)(\textfrak{b}'\textfrak{b})^{-1}$ will live in $\textfrak{B}_{\rho^{(i)}}\textfrak{U}$. Since $[\mathbb{Z}\rho^{(i_t)} + \mathbb{Z}][\mathbb{Z}\rho^{(i)} + \mathbb{Z}] = [\mathbb{Z}\rho^{(j)} + \mathbb{Z}]$, we have the product $(r_2\rho^{(i_t)} - r_1)(s_2\rho^{(i)} - s_1) = (u_2\rho^{(j)} - u_1)$ for some

\begin{align}
    \mathbf{u} &= \begin{bmatrix}
        u_1 \\
        u_2
    \end{bmatrix} \in \mathbb{Q}^2
\end{align} 

\noindent
We conclude that

$$h(\exp(n Z'_{d \infty_2}(0,\textfrak{U}))) = \exp(n Z'_{d\infty_2}(0,\textfrak{B}_{\rho^{(i)}}\textfrak{U})) = (\psi^{-2}\chi_{\mathbf{u}}^{-1})(A_t)\shin^{\mathbf{u}}_A(\rho^{(j)})^2$$

\noindent
implying:

\begin{theorem}
Assume the order 1 abelian Stark conjecture for real quadratic fields. Then using the notation defined above we have
$$h(\phi_\mathbf{p}(t)\shin_A^{d^{-1}\mathbf{p}}(\rho_t)) = \pm \phi_{d\mathbf{u}}(t)\shin_A^\mathbf{u}(\rho^{(j)})$$
whenever $h$ corresponds to a non-principal ideal in $\text{Cl}_{\overline{d} \infty_2}(\mathcal{O}_f)$.
\end{theorem}

\begin{conjecture}
We conjecture that the sign in both cases discussed in this note is equal to $h(\sqrt{d+1})/\sqrt{d+1}$ so that using the notation defined above we have
\[
h(\phi_\mathbf{p}(t)\shin_A^{d^{-1}\mathbf{p}}(\rho_t)) = \frac{h(\sqrt{d+1})}{\sqrt{d+1}}\phi_{M_b\mathbf{p}}(t)\shin_A^{d^{-1}M_b\mathbf{p}}(\rho_t)
\]
whenever $h$ corresponds to a principal ideal and 
\[
h(\phi_\mathbf{p}(t)\shin_A^{d^{-1}\mathbf{p}}(\rho_t)) = \frac{h(\sqrt{d+1})}{\sqrt{d+1}} \phi_{d\mathbf{u}}(t)\shin_A^\mathbf{u}(\rho^{(j)})
\]
whenever $h$ corresponds to a non-principal ideal.
\end{conjecture}

\noindent
Numerical evidence for this conjecture is provided by the algorithm which is described below. An implementation of this algorithm can be found \href{https://github.com/jddicapua1023/GaloisZauner}{here}. An AI tool was used solely to implement this code. The mathematics (including the algorithm) and this paper's text are my own. This implementation is verified for $d=9,15$ and $d=21$ through $d=70$ as of 11/3/2026.

\section{Description of the algorithm}
\label{sec4}

At this point we give a description of our algorithm which takes an admissible tuple attached to a dimension $d$ and outputs the Galois action of the ray class group $\text{Cl}_{\overline{d}\infty_2}(\mathcal{O}_f)$ on the ghost SIC attached to the admissible tuple via Zauner.jl \cite{ApplebyFlammiaKopp}. The algorithm starts by generating tables of representatives for $\text{Cl}_{\overline{d}\infty_2}(\mathcal{O}_f)$ and $\overline{\text{Clm}}^\flat_{\overline{d}\infty_2}(\mathcal{O}_f)$. For each entry in these tables, the code records the pair $(\rho,\mathbf{r})$ furnished by Theorem 2.19 \cite{ApplebyFlammiaKopp} and the corresponding \say{$\shin$ value} at level $\overline{d}$. It then takes the ghost SIC $\tilde{\Pi}$ attached to the admissible tuple and constructed by Zauner.jl and a representative $J$ for a class of $\text{Cl}_{\overline{d}\infty_2}(\mathcal{O}_f)$. For each ghost overlap $\tilde{\mu}_\mathbf{p}$ of $\tilde{\Pi}$, the code identifies the class $\mathfrak{U} \in \overline{\text{Clm}}^\flat_{\overline{d}\infty_2}(\mathcal{O}_f)$ attached to $\tilde{\mu}_\mathbf{p}$ by Theorem 2.19, namely the class of the ideal generated by the $\Sigma$-positive lift of $p_2\rho_t-p_1$. This identification is made by exact membership in the monoid table, and the code checks that the \say{$\shin$ value} tabulated for $\mathfrak{U}$ agrees with the one computed at $\mathbf{p}$. It then computes the (conjectured) Galois action of $J$ on $\tilde{\mu}_\mathbf{p}$ via the method described in \hyperref[sec3]{Section 3} of this paper. Suppose $\text{Art}([J]) = h \in \text{Gal}(H^{\mathcal{O}_f}_{\overline{d}\infty_2}/K)$. Then the code reconstructs $h(\tilde{\Pi})$ from the table of values $h(\tilde{\mu})$.

\end{document}